# A Group-Permutation Algorithm to Solve the Generalized SUDOKU

Florentin Smarandache
University of New Mexico
Gallup Campus, USA

<u>Sudoku</u> is a game with numbers, formed by a square with the side of 9, and on each row and column are placed the digits 1, 2, 3, 4, 5, 6, 7, 8, 9, written only one time; the square is subdivided in 9 smaller squares with the side of $3 \times 3$, which, also, must satisfy the same condition, i.e. each square to contain all digits from 1 to 9 written only once.

The Japanese company Nikoli has popularized this game in 1986 under the name of *sudoku*, meaning "single number".

Sudoku can be generalized to squares whose dimensions are $n^2 \times n^2$, where $n \geq 2$, using various symbols (numbers, letters, mathematical symbols, etc.), written just one time on each row and on each column; and the large square is divided into $n^2$ small squares with the side $n \times n$ and each will contain all $n^2$ symbols written only once.

An <u>elementary solution</u> of one of these <u>generalized Sudokus</u>, with elements (symbols) from the set

$$S = \{s_1, s_2, ..., s_n, s_{n+1}, ..., s_{2n}, ..., s_{n^\wedge 2}\}$$

(supposing that their placement represents the relation of total order on the set of elements $S$), is:

Row 1: all elements in ascending order

$$s_1, s_2, ..., s_n, s_{n+1}, ..., s_{2n}, ..., s_{n^\wedge 2}$$

On the next rows we will use circular permutations, considering groups of $n$ elements from the first row as follows:

Row 2:
$$s_{n+1}, s_{n+2}, ..., s_{2n}; s_{2n+1}, ..., s_{3n}; ...; s_{n^\wedge 2}; s_1, s_2, ..., s_n$$

Row 3:
$$s_{2n+1}, ..., s_{3n}; ...; s_{n^\wedge 2}; s_1, s_2, ..., s_n; s_{n+1}, s_{n+2}, ..., s_{2n}$$

...........................................................................................

Row $n$:
$$s_{n^\wedge 2 - n+1}, ..., s_{n^\wedge 2}; s_1, ..., s_n, s_{n+1}; s_{n+2}, ..., s_{2n}; ...; s_{3n}; ...; s_{n^\wedge 2 - n}.$$

Now we start permutations of the elements of row $n+1$ considering again groups of $n$ elements.

Row $n+1$:
$$s_2, ..., s_n, s_{n+1}; s_{n+2}, ..., s_{2n}, s_{2n+1}; s_{n^\wedge 2 - n+2}, ..., s_{n^\wedge 2}, s_1$$

Row $n+2$:
$$s_{n+2}, ..., s_{2n}, s_{2n+1}; s_{n^\wedge 2 - n+2}, ..., s_{n^\wedge 2}, s_1; s_2, ..., s_n, s_{n+1}$$

……………………………………………………………………….
Row $2n$:

$S_{n^2-n+2}$ ..., $S_{n^2}$, $S_1$; $S_2$, ..., $S_n$, $S_{n+1}$; $S_{n+2}$, ..., $S_{2n}$, $S_{2n+1}$

Row $2n+1$:

$S_3$, ..., $S_{n+2}$; $S_{n+3}$, ..., $S_{2n+2}$; $S_{n^2+3}$, ..., $S_{n^2}$, $S_1$, $S_2$

and so on.

Replacing the set $S$ by any permutation of its symbols, which we'll note by $S'$, and applying the same procedure as above, we will obtain a new solution.

The classical Sudoku is obtained for $n = 3$.

Below is an example of this group-permutation algorithm for the classical case:

| 1 | 2 | 3 | 4 | 5 | 6 | 7 | 8 | 9 |
|---|---|---|---|---|---|---|---|---|
| 4 | 5 | 6 | 7 | 8 | 9 | 1 | 2 | 3 |
| 7 | 8 | 9 | 1 | 2 | 3 | 4 | 5 | 6 |
| 2 | 3 | 4 | 5 | 6 | 7 | 8 | 9 | 1 |
| 5 | 6 | 7 | 8 | 9 | 1 | 2 | 3 | 4 |
| 8 | 9 | 1 | 2 | 3 | 4 | 5 | 6 | 7 |
| 3 | 4 | 5 | 6 | 7 | 8 | 9 | 1 | 2 |
| 6 | 7 | 8 | 9 | 1 | 2 | 3 | 4 | 5 |
| 9 | 1 | 2 | 3 | 4 | 5 | 6 | 7 | 8 |

For a $4^2 \times 4^2$ **square** we use the following 16 symbols:
{A, B, C, D, E, F, G, H, I, J, K, L, M, N, O, P}
and use the same group-permutation algorithm to solve this Sudoku.

From one solution to the generalized Sudoku we can get more solutions by simply doing permutations of columns or/and of rows of the first solution.

| A | B | C | D | E | F | G | H | I | J | K | L | M | N | O | P |
|---|---|---|---|---|---|---|---|---|---|---|---|---|---|---|---|
| E | F | G | H | I | J | K | L | M | N | O | P | A | B | C | D |
| I | J | K | L | M | N | O | P | A | B | C | D | E | F | G | H |
| M | N | O | P | A | B | C | D | E | F | G | H | I | J | K | L |
| B | C | D | E | F | G | H | I | J | K | L | M | N | O | P | A |
| F | G | H | I | J | K | L | M | N | O | P | A | B | C | D | E |
| J | K | L | M | N | O | P | A | B | C | D | E | F | G | H | I |
| N | O | P | A | B | C | D | E | F | G | H | I | J | K | L | M |
| C | D | E | F | G | H | I | J | K | L | M | N | O | P | A | B |
| G | H | I | J | K | L | M | N | O | P | A | B | C | D | E | F |
| K | L | M | N | O | P | A | B | C | D | E | F | G | H | I | J |
| O | P | A | B | C | D | E | F | G | H | I | J | K | L | M | N |
| D | E | F | G | H | I | J | K | L | M | N | O | P | A | B | C |
| H | I | J | K | L | M | N | O | P | A | B | C | D | E | F | G |
| L | M | N | O | P | A | B | C | D | E | F | G | H | I | J | K |
| P | A | B | C | D | E | F | G | H | I | J | K | L | M | N | O |